\documentclass[11pt]{article}
\usepackage{hyperref}
\usepackage{amssymb}
\usepackage{amsmath}
\usepackage{amsthm}
\usepackage{graphicx}
\usepackage{subfigure,url,breakurl}
\usepackage{bridges}

\newtheorem{theorem}{Theorem}

\newtheorem{lemma}{Lemma}

\newtheorem{problem}{Problem}

\newcommand{\cref}[1]{Corollary \ref{cor:#1}}

\newcommand{\fref}[1]{Fig. \ref{fig:#1}}

\newcommand{\probref}[1]{Problem \ref{prob:#1}}

\renewcommand{\subsection}{\paragraph}

\title{Tatami Maker: A combinatorially rich mechanical game board}
\author{Alejandro Erickson\\
  Department of Computer Science, University of Victoria, V8W 3P6,
  Canada}
\date{}

\begin{document}



  \maketitle



  \thispagestyle{empty}
  \begin{abstract}
    Japanese tatami mats are often arranged so that no four mats meet.
    This local restriction imposes a rich combinatorial structure when
    applied to monomino-domino coverings of rectilinear grids.  We
    describe a modular, mechanical game board, prototyped with a
    desktop 3D printer, that enforces this restriction, and transforms
    tatami pen-and-paper puzzles into interactive sculptures. We
    review some recent mathematical discoveries on tatami coverings
    and present five new combinatorial games implemented on the game
    board.
\end{abstract}

\section*{Introduction}
\label{sec:intro}

Tatami mats are a common floor furnishing, originating in aristocratic
Japan, during the Heian period (794-1185).
These thick mats, once hand-made with a rice straw core and a soft,
woven rush straw exterior, are now machine-produced in a variety of
materials, and are available in mass-market stores. They are so
integral to Japanese culture, that a standard sized mat is the unit of
measurement in many architectural applications (see \cite{Hiroshi1983}).

In our mathematical musings we depart considerably from traditional
layouts, but we retain two essential items.  The first of these is
aspect ratio; a full mat is a $1\times 2$ domino, and a half mat is
$1\times 1$ monomino.  The second is the 17th century rule for
creating auspicious arrangements; \emph{no four mats may meet}.


Counting domino coverings is a classic area of enumerative
combinatorics and theoretical computer science, but little attention
has been paid to problems where the local interactions of the dominoes
are resctricted in some fashion.  The tatami restriction is perhaps
the most natural of these, and it imposes a visually appealing
structure with nice combinatorial properties (see \fref{tatami}).  As
a result, it is the subject of an exercise in Volume 4 of ``The Art of
Computer Programming'' (\cite{Knuth2011}), where Knuth reprints a
diagram of a 17th century Japanese mathematician, and recently the
tatami restriction has been studied in several research papers (see
\cite{ALHAZOVMORITAIWAMOTO2010,EricksonRuskeySchurch2011,EricksonSchurch2012,RuskeyWoodcock2009,EricksonRuskey2013a,EricksonRuskey2013}).

\begin{figure}[h]
  \centering
  \includegraphics[width=\textwidth]{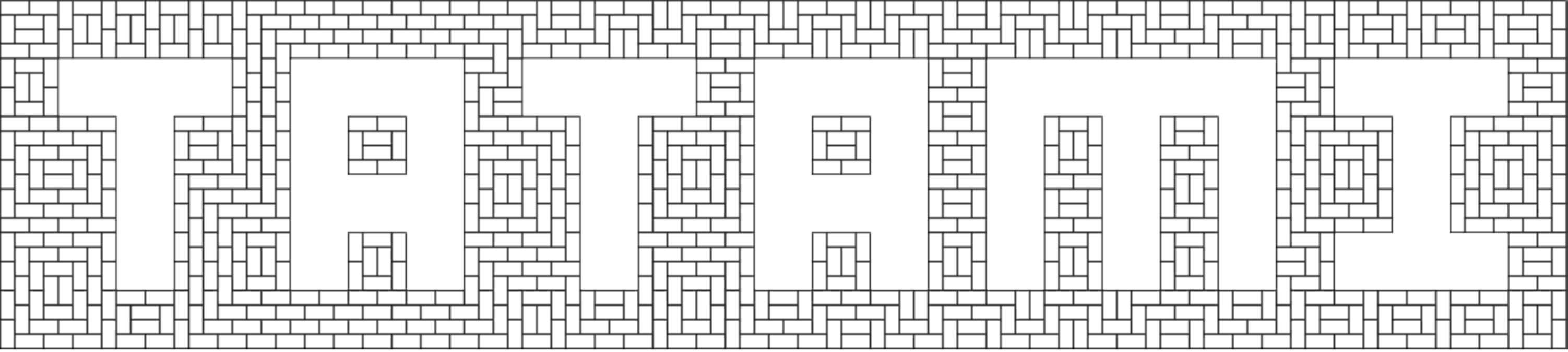}
  \caption{A tatami domino covering, generated in software using a
    logical satisfiability formula.}
  \label{fig:tatami}
\end{figure}

Mathematics in art is an ideal pedagogic formula; the math challenges
the logical capabilities of the participant, whereas the art compels
her to understand it.  The combination is even more effective for an
interactive art installation that changes and increases gradually in
sophistication as the participant furthers her engagement with it.  An
example of this is an installation called \emph{Boundary Functions},
by Scott Snibbe (see~\cite{Snibbe1998}).  Participants
on the floor of the installation become points in a Voronoi diagram
that is projected onto the floor from the ceiling.  That is, lines are
projected to form a 2-dimensional cell around each participant.  As
more participants enter the installation and move around, the
installation reveals more about the nature of Voronoi diagrams.
Boundary Functions may or may not have been intended as a teaching
tool, but learning would appear to be an inevitable consequence.


Snibbe's installation is inviting, aesthetically appealing, and self
explanatory.  I struggle to envision a Voronoi-themed pen-and-paper
puzzle that qualifies as art, precisely because it would be missing
these things.  The same challenge presents itself for tatami
coverings.  We describe four tatami puzzle games which have the
potential to become art by means of a kinetic sculpture.

On paper, a tatami puzzle is comprised of a grid and some
instructions. The player must internalize the tatami restriction, and
apply it every time she draws a tile.  The structure of a tatami
covering is not obvious to the uninformed, and applying the local
restriction can overwhelm the player's efforts to solve the puzzle.
Furthermore, the appearance of the pen-and-paper puzzle depends on the
player's artistic ability and care, which can negate the natural
beauty of tatami coverings.

A mechanical game board that enforces the tatami restriction
eliminates the preoccupation with ``no four tiles meet''.  The
uninformed player engages with the installation by trial and error and
the sculpture teaches the tatami restriction incrementally.  At the
end she is rewarded with a completed tatami covering.

We describe the structure of tatami coverings and a prototypical
kinetic sculpture called Tatami Maker, which realizes all of the
combinatorial properties of tatami coverings.  We introduce four
puzzle games that are playable on Tatami Maker, along with some
related mathematical results.

\section*{Structure}

We introduce the tatami structure with an interactive demonstration.
Reader, take up your pencil and complete the partial covering in
\fref{comprehensiveSource}.  Once again, no four mats may meet at any
point; alternatively, every intersection of the grid will contact a
broad edge of a domino.  You should find the completion is unique and
equivalent to \fref{comprehensiveTiling}.

Completing the exercise will inevitably bring about the discovery of
the \emph{ray}, which occurs wherever a vertical domino shares an edge
with a horizontal domino.  The tatami restriction forces the pattern
to repeat itself until it reaches the boundary of the grid.

\begin{figure}[h]
  \centering
  %
  \subfigure[ ]{%
    \includegraphics[width=0.37954\textwidth]{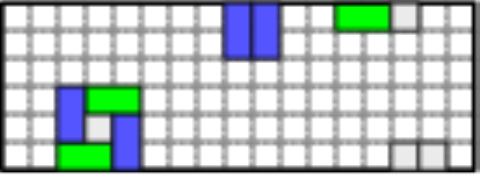}%
    \label{fig:comprehensiveSource}}%
  \hspace{0.01\textwidth}%
  \subfigure[ ]{%
    \includegraphics[width=0.37954\textwidth]{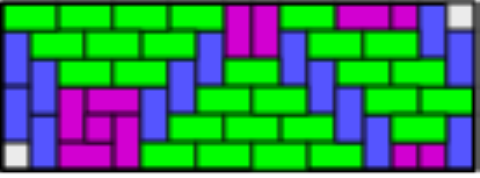}%
    \label{fig:comprehensiveTiling}}%
  \hspace{0.01\textwidth}%
    \subfigure[ ]{%
    \includegraphics[width=0.20093\textwidth]{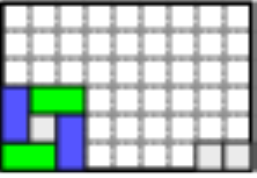}%
    \label{fig:incompatibleSources}}%
  \caption{Attempt to complete the partial tatami coverings by placing
    tiles which are forced by the tatami restriction.  \emph{(a)}, A partial
    covering and, \emph{(b)}, its unique completion.  \emph{(b)}, This covering
    contains every possible type of feature.  \emph{(c)} This partial
    covering cannot be completed.}
  \label{fig:comprehensive}
\end{figure}

The tiles in the partial covering in \fref{incompatibleSources} are
incompatible because the tatami restriction forces the propagation of
rays that cross paths.  I invite the reader to check.

Tatami coverings have the remarkable property that a simple local
restriction imposes an aesthetically appealing, combinatorial
structure.  We discovered the ray in \fref{comprehensive}, and a case
analysis of how rays begin reveals that in a rectangular covering,
there are only four distinct ray-producing \emph{features}, up to
reflection and rotation.  These four features and their rays are shown
in \fref{features} (see \cite{EricksonRuskeySchurch2011}). A brick
laying pattern, called \emph{bond}, fills the areas between the
features and rays.  There are two types of horizontal bond and two
types of vertical bond.

\begin{figure}[h]
  \centering
  \subfigure[loner]{%
    \includegraphics[scale=0.5]{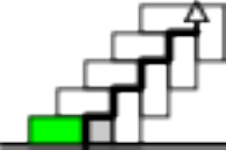}%
    \label{fig:loner}}%
  \hspace{0.5in}%
  \subfigure[vee]{%
    \includegraphics[scale=0.5]{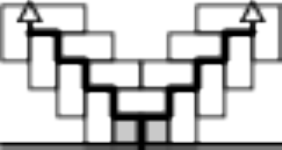}%
    \label{fig:vee}}%
  \hspace{0.5in}
  \subfigure[bidimers]{%
  \includegraphics[scale=0.5]{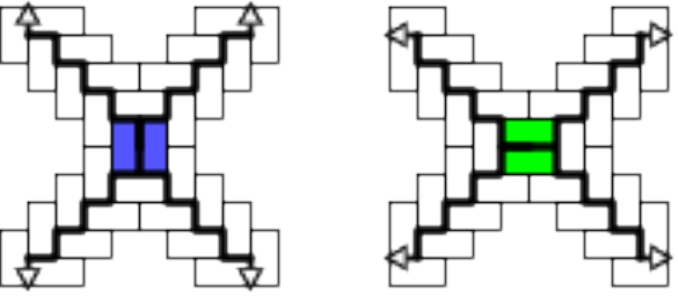}%
    \label{fig:bidimer}}%
  \hspace{0.5in}%
  \centering \subfigure[vortices]{%
    \includegraphics[scale=0.5]{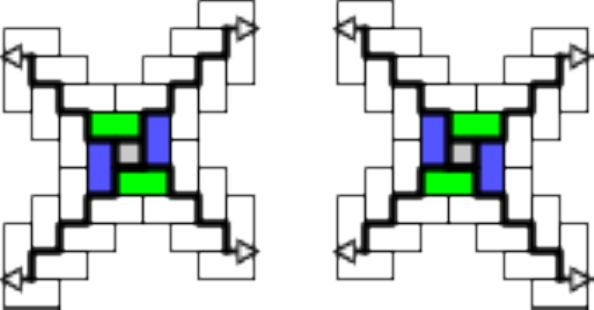}%
    \label{fig:vortex}}%

  \caption{Tatami coverings of rectangular grids are comprised of
    these four types of feature and the bond pattern.}
  \label{fig:features}
\end{figure}

Likely the following theorem holds for general rectilinear regions,
but it has only been proved for rectangular grids.

\begin{lemma}[\cite{EricksonRuskeySchurch2011}]
  \label{lem:boundary}
  Let $G$ be an $r\times c$ grid, with $r<c$.  A tatami covering of $G$
  is uniquely determined by the tiles touching its boundary.
\end{lemma}

In terms of the area of the region, tatami coverings are simpler than
general monomino-domino coverings by an entire order of complexity.
The result is that a randomly chosen tatami covering has an inevitable
simplicity that is beautiful, not only in its construction, but also
in its appearance.  In spite of this, the structure described here is
not obvious to the uninformed, and building an intuition for it can
take considerable effort.

\section*{Tatami Maker: a combinatorially rich mechanical game board}
The Tatami Maker game board forms a rectilinear grid that enforces the
tatami restriction when it is covered by the accompanying tile pieces.
Arbitrary rectilinear grids can be created by placing Tatami Maker's
modules alongside each other.  A simple mechanism is embedded at each
grid line intersection, which obstructs the placement of a tile if the
other three incident grid squares are covered.

\begin{figure}[h]
  \centering
  \subfigure[Tatami Maker modules are placed alongside each other to
  create larger grids.]{%
    \includegraphics[width=0.35\textwidth]{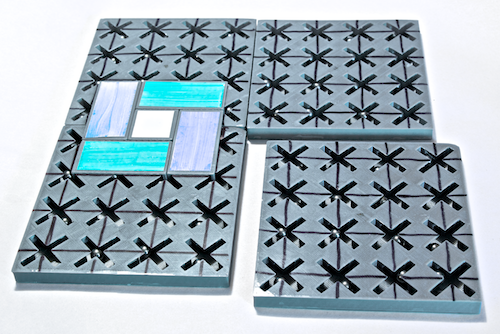}%
    \label{fig:modules}}%
  \hspace{0.5in}%
  \subfigure[The inside of Tatami Maker's intersection mechanism.
  Each arm of the \texttt{X} shape extends to a different grid
  square.]{%
    \includegraphics[width=0.30\textwidth]{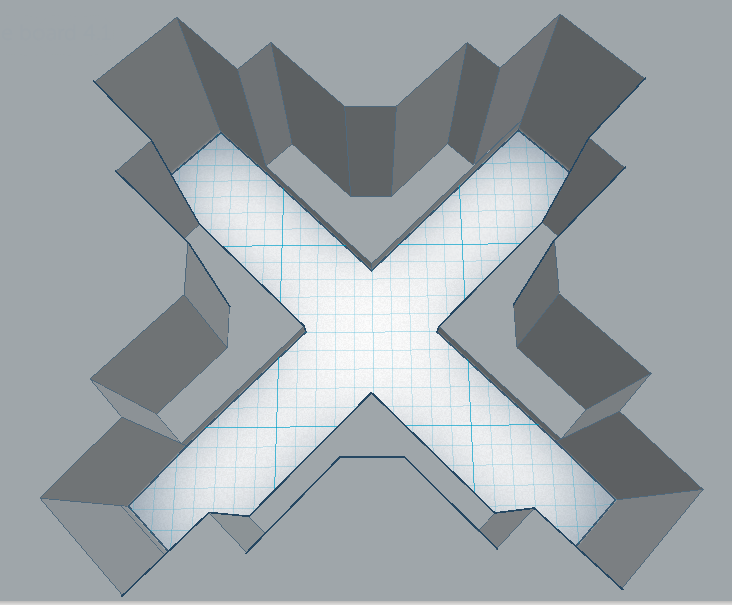}%
    \label{fig:mechanism}}%
  \caption{Tatami Maker modules and mechanism.}%
  \label{fig:mechanismandmodules}
\end{figure}

A tile piece covers one or two grid squares, and the underside of each
of its four corners has the specially shaped foot shown in
\fref{tilepiece}.  The feet interact with the mechanism in the game
board by pushing an obstructing ball onto an unoccupied grid square,
as well as guiding the tile correctly onto the grid.

Each mechanism occupies one grid intersection, which consists of one
quadrant from each of four incident grid squares.  A cavity in the
game board contains a ball which may travel freely to any unoccupied
quadrant.  If a quadrant is occupied by a tile's corner, then one of
tile's feet occupies the part of the cavity that is otherwise
available to the ball.  When three of the quadrants are occupied by
tile corners, the ball is forced to occupy the remaining quadrant,
thereby preventing a tile corner from being placed here.

A minimal game board module is one grid intersection.  As a result,
the dimensions of a module are given in terms of its intersections
rather than its grid squares.  Our prototype's modules are $4\times
4$, but they need not be square, or even rectangular.

\begin{figure}[h]
  \centering
  \subfigure[ ]{%
    \includegraphics[width=0.30\textwidth]{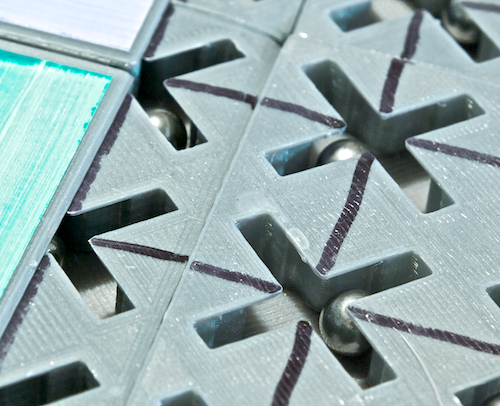}%
    \label{fig:tatamimaker}}%
  \hspace{0.5in}%
  \subfigure[ ]{%
    \includegraphics[width=0.40\textwidth]{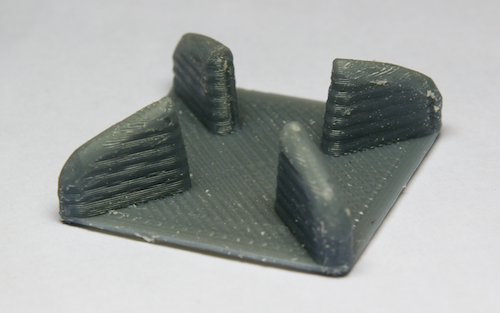}%
    \label{fig:tilepiece}}%
  \caption{\emph{(a)}, Grid intersections and several Tatami Maker modules.
    \emph{(b)} The feet of a tile piece.}
  \label{fig:makerandpiece}
\end{figure}

The main design challenge is ensuring that the ball can be pushed by
an incoming foot, unobstructed, to an available quadrant.  We label
the quadrants Q1, Q2, Q3, and Q4, in counterclockwise order.  A
critical case occurs when Q2 and Q4 are occupied by tile corners (and
feet), and a tile foot placed in Q1 must push the ball to Q3.
Intuitively, the ball must disturb the feet in Q2 and Q4 on the way to
Q3, otherwise the midpoint of the intersection would accomodate the
ball when all four quadrants are occupied.  The mechanism is designed,
therefore, so that the ball lifts the tiles in Q2 and Q4 as it passes
to Q3, and so that it will not become stuck against another part of
the mechanism before it arrives in Q3.

The availability of desktop 3D printers has lowered the cost of
creating prototypes sufficiently that rough estimation, and iterative
trial and error was the most economical way of solving these design
challenges.  Tatami Maker was prototyped with a Solidoodle 2 printer;
a 3D CNC machine, that extrudes a filament of hot ABS plastic into a
$15 \times 15 \times 15$ cm print area to create real life versions of
virtual 3D models.  With each iteration of the prototype, changes to
the shape of the cavity and the feet of the tile pieces were made to
tease out the required behaviour of the mechanism.

Tatami Maker was tested informally at Science World in Vancouver, and
it mesmerized many passing young children (see \fref{scienceworld}).

\begin{figure}[h]
  \centering
  \subfigure[ ]{%
    \includegraphics[width=0.35\textwidth]{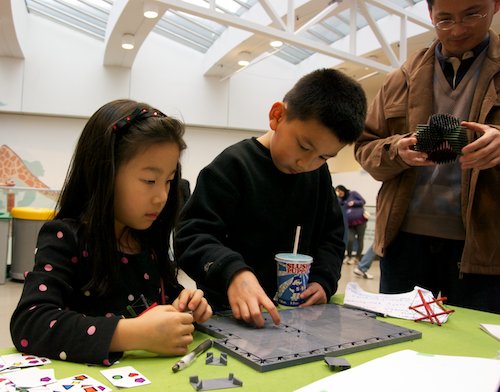}%
    \label{fig:scienceworld1}}%
  \hspace{0.5in}%
  \subfigure[ ]{%
    \includegraphics[width=0.35\textwidth]{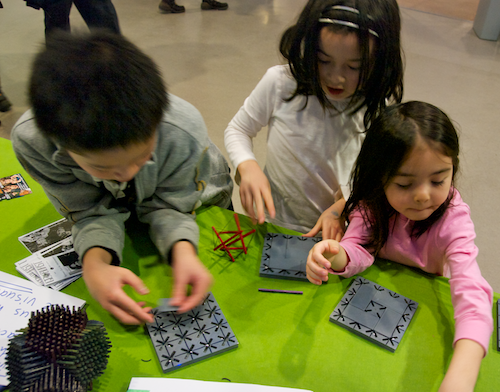}%
    \label{fig:scienceworld2}}%
  \caption{Visitors to Science World in Vancouver, Canada, test Tatami
    Maker by playing Oku.}
  \label{fig:scienceworld}
\end{figure}

\section*{Five tatami puzzles}
\label{sec:puzzles}

We describe how to set up five original combinatorial puzzle games for
Tatami Maker, and highlight some related mathematical discoveries.
The first four are for a single player, and the last one is for two
players.

\subsection*{Oku}
From the Japanese word for ``put'', the most straightforward puzzle
requires that the player cover a rectilinear region with tiles.  The
informed player can do this very easily by using a bond pattern, and
placing monominoes wherever dominoes do not fit, but otherwise the
pitfalls arising from the structure of tatami coverings are numerous
(see \fref{oku}), and Oku is particularly instructive in this regard.

\begin{figure}[h]
  \centering
  \includegraphics[width=0.5\textwidth]{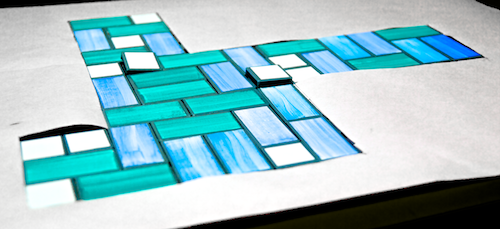}
  \caption{Tatami Maker modules lie under an overlay.  The two raised
    monominoes are obstructed by Tatami Maker's mechanism.}
  \label{fig:oku}
\end{figure}
Some of the nicest combinatorial results for Oku are on square grids.
The numbers of coverings for various parameters have simple
representations, and divide up into natural equivalence classes (see
\cite{EricksonSchurch2012}).

\begin{theorem}[\cite{EricksonRuskeySchurch2011,
    EricksonSchurch2012}] \label{thm:square} If $n$ and $m$ have the
  same parity, then the number of tatami coverings of the $n\times n$
  grid with $m$ monominoes is
  \begin{itemize}
  \item $m2^{m} + (m+1)2^{m+1}$ when $m<n$,
  \item $n2^{n-1}$ when $m=n$; and
  \item $0$ when $m>n$.
  \end{itemize}
\end{theorem}



We can also generate certain subsets of these coverings, for example
those with $m=n$ and a given number of vertical dominoes (see
\fref{n8k7}).


\begin{figure}[h]
  \centering
\includegraphics[width=0.5\textwidth]{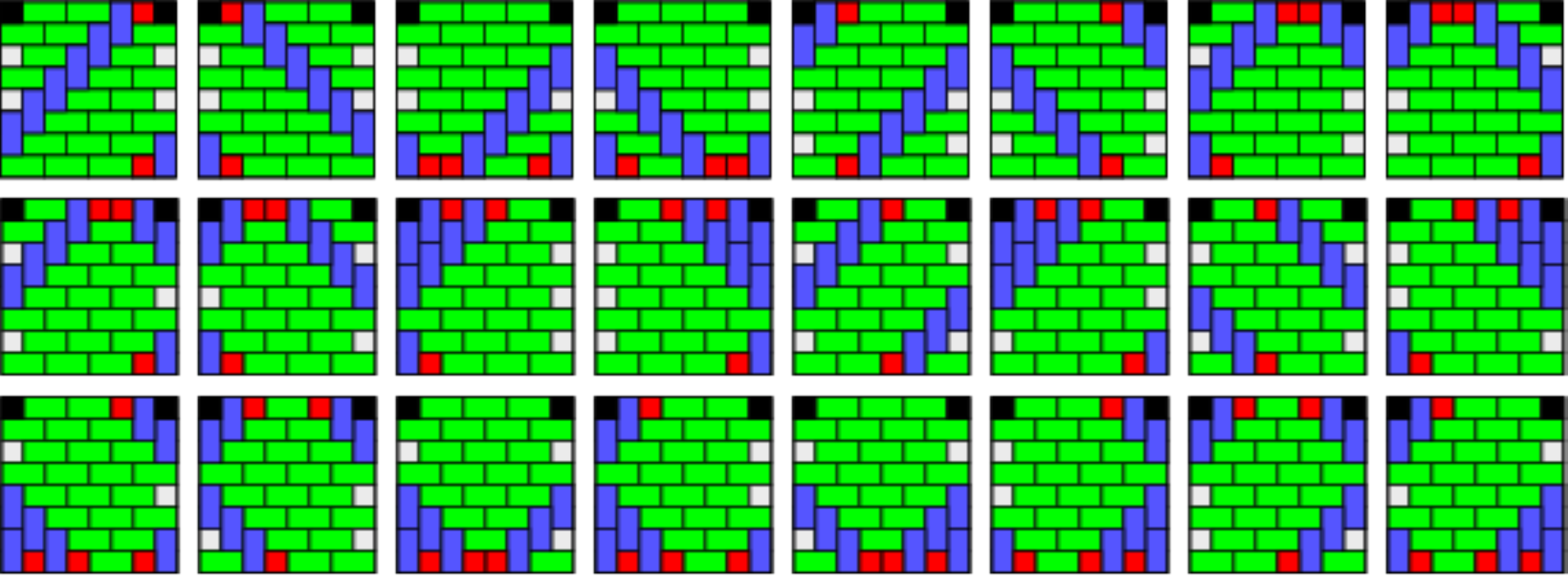}  
\caption{All $8\times 8$ coverings with $7$ vertical dominoes.}
  \label{fig:n8k7}
\end{figure}

\subsection*{Tomoku}
\label{sec:tomoku}
Due to Mart\'in Matamala (private communication), Tomoku is a
contraction of the words tomography, and the above game, Oku.  This
puzzle is played on a rectangular board, and the player is given the
row and column projection of a set of solutions, one of which the
player must find to complete the puzzle.

\begin{figure}[h]
  \centering
  \subfigure[]{%
    \includegraphics[width=0.475\textwidth]{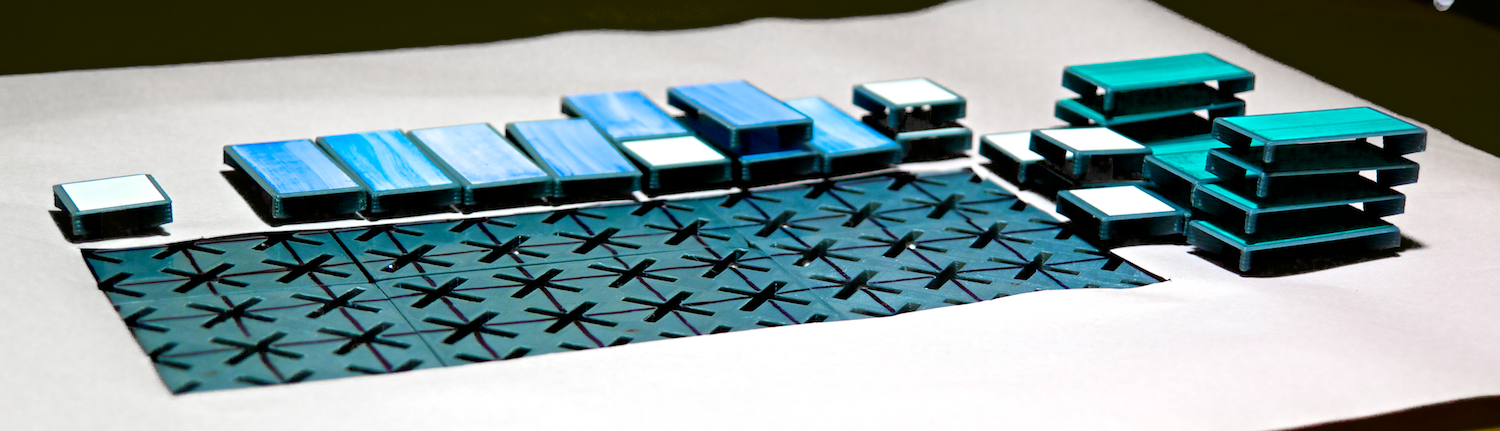}%
    \label{fig:tomokusetup}}%
  \hspace{0.2in}%
  \subfigure[]{%
    \includegraphics[width=0.475\textwidth]{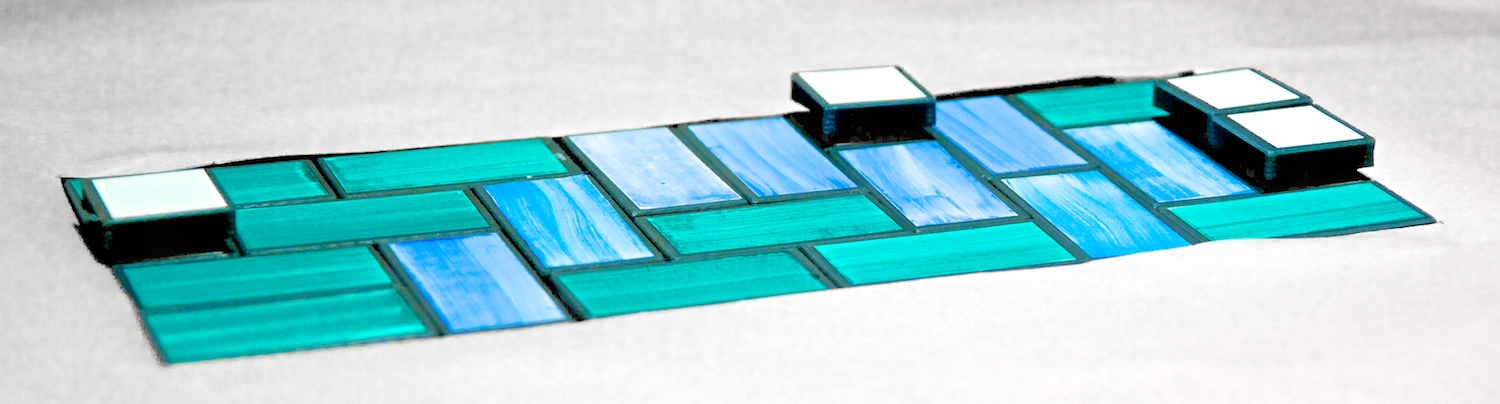}%
    \label{fig:tomokusolved}}%
  \caption{\emph{(a)}, Typical setup for an instance of Tomoku.  Monominoes are
    double stacked in solution, \emph{(b)} because they appear in both a
    column and a row.}
\label{fig:tomoku}
\end{figure}

Tomoku is defined as a decision problem, where one must determine
whether or not a solution exists.  In practice, it is more
interesting to reconstruct a covering when the answer is yes, than it
is to solve a decision problem.

\begin{problem}
  \label{prob:tomoku}
  Given an $r\times c$ grid, and a triple of integers, $(v,h,m)$, for
  each row and each column, denoting the number of grid squares
  covered by vertical dominoes, horizontal dominoes, and monominoes,
  determine whether or not a covering exists with these row and column
  projections.
\end{problem}

It is not known whether a polynomial algorithm exists to solve an
instance of \probref{tomoku}, but there are pairs of coverings with
the same row and column projection.  For example, an $n\times n$
covering with a central clockwise vortex gives the same row and column
projections as a counterclockwise vortex.

An instance of Tomoku is set up by arranging Tatami Maker modules
under an overlay which is printed with instructions and masks the
modules to expose only the required grid.  The instructions indicate
how to set up the modules and overlay, as well as the details of one
or more instances of Tomoku (see \fref{tomokuInstructions}).

\begin{figure}[h]
  \centering
  \includegraphics[width=\textwidth]{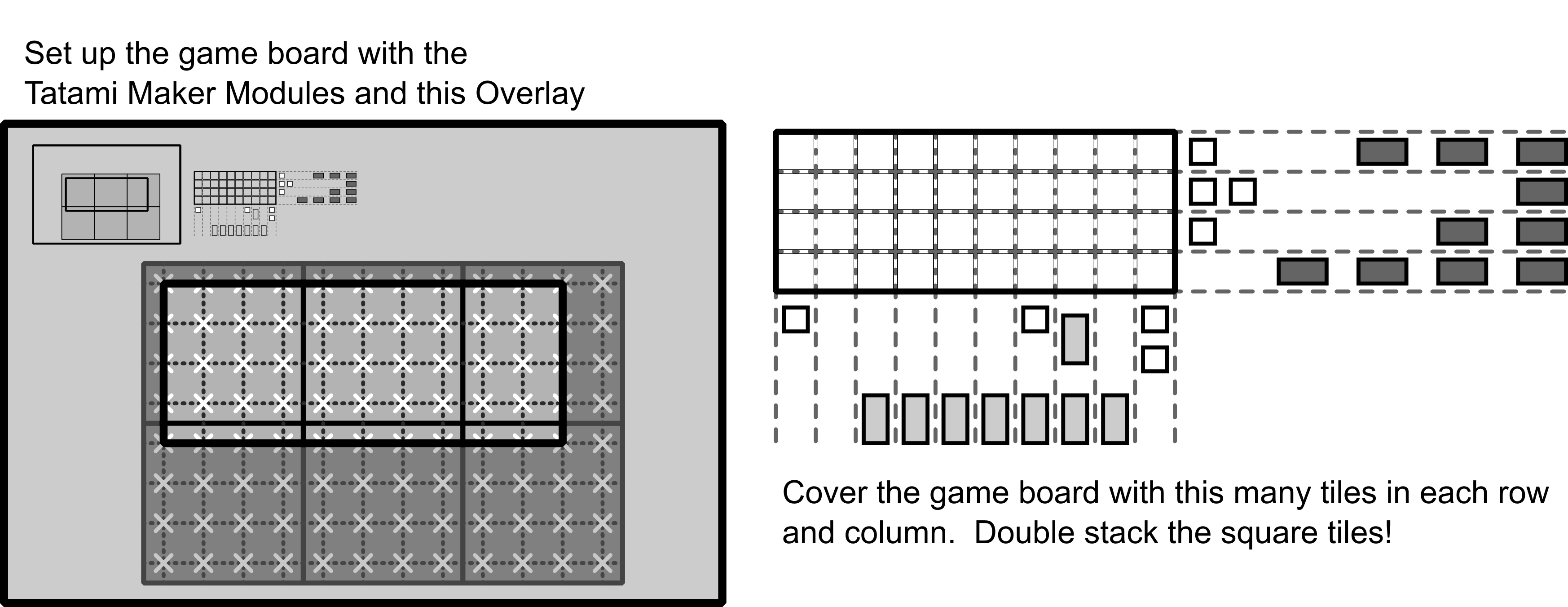}
  \caption{Instructions printed on the Tomoku overlay for the instance
    in \fref{tomokusetup}.  This particular one can be completed
    without backtracking.}
  \label{fig:tomokuInstructions}
\end{figure}

Tomoku may be played with pencil and paper, and two instances are
reprinted in \fref{tomokuExamples} from \cite{Erickson2012}. Considerable
efficiency can be acheived by using short line segments to represent
dominoes, and dots to represent monominoes.

\begin{figure}[h]
  \centering
  \subfigure[ ]{%
    \includegraphics[width=0.475\textwidth]{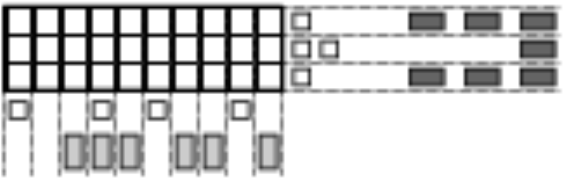}%
    \label{fig:r3c10p1}}%
  \hspace{0.2in}%
  \subfigure[ ]{%
    \includegraphics[width=0.475\textwidth]{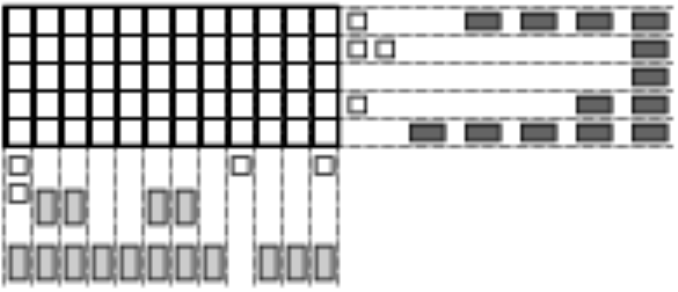}%
    \label{fig:r5c12p1}}%
  \caption{Instances of Tomoku, reprinted from
    \cite{Erickson2012}. The $5\times 12$ puzzle is quite
    challenging.}
  \label{fig:tomokuExamples}
\end{figure}

\subsection*{The Lazy Paver}
\label{sec:lazy}
The Lazy Paver is tasked with tatami covering a rectilinear driveway
with pavers.  She would typically produce a $1\times 1$ paver by
cutting a $1\times 2$ paver in half, so instead, she will avoid the
extra work by covering the driveway with domino pavers.

\begin{figure}[h]
  \centering
  \subfigure[ ]{%
  \includegraphics[width=0.65\textwidth]{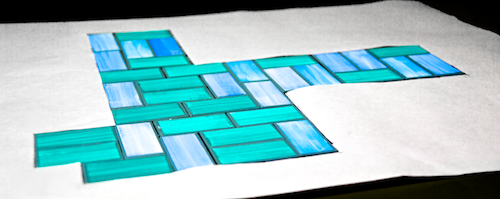}%
    \label{fig:lazy}}%
  \hspace{0.2in}%
  \subfigure[ ]{%
    \includegraphics[width=0.3\textwidth]{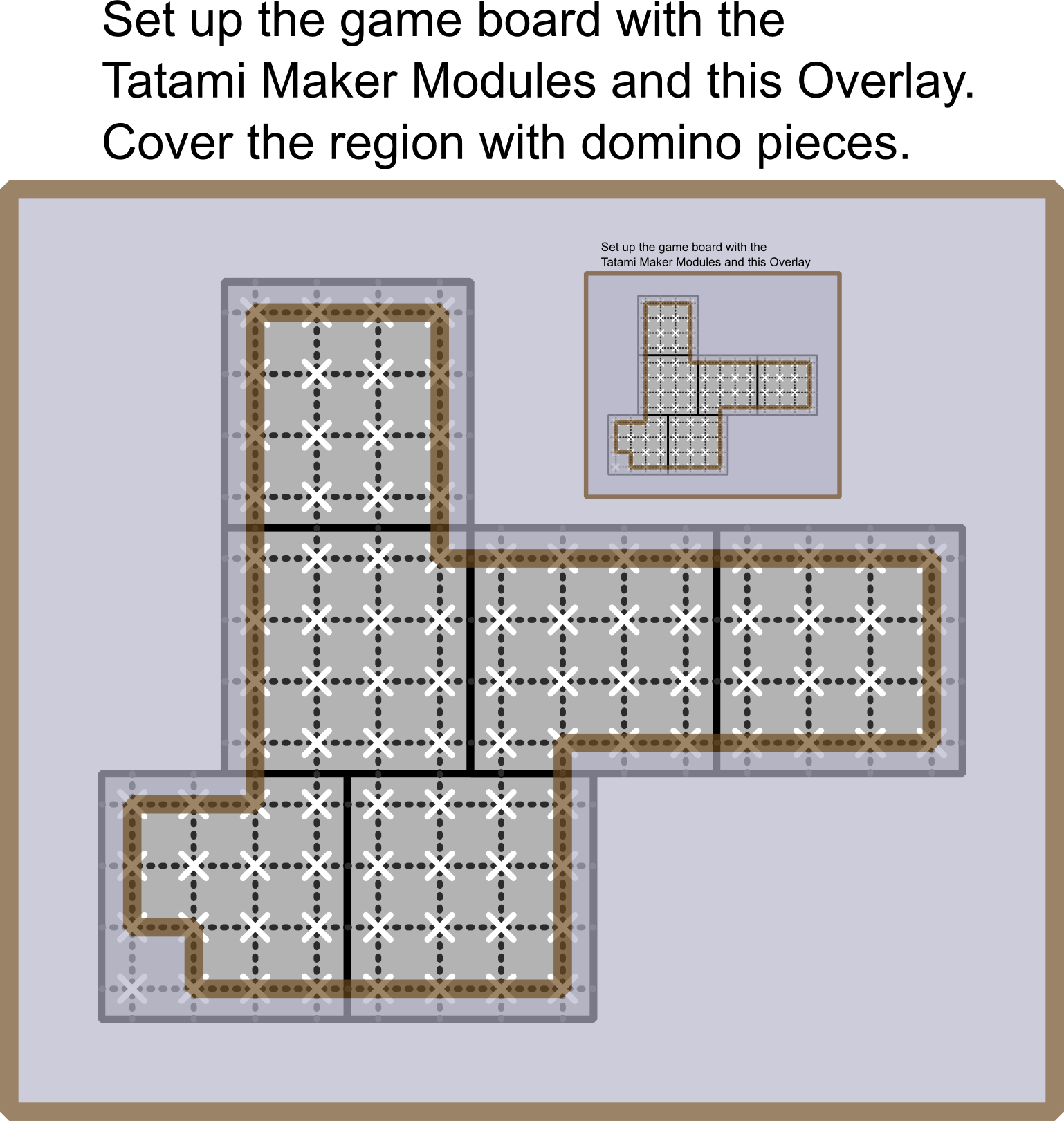}%
    \label{fig:lazyPaverInstructions}}%
  \caption{Lazy Paver instance and instructions printed on an overlay.}
  \label{fig:lazyPaver}
\end{figure}

The Lazy Paver is also defined as a decision problem.

\begin{problem}[\cite{EricksonRuskeySchurch2011}]
  \label{prob:lazy}
  Determine whether or not a given rectilinear region can be tatami
  covered with dominoes.
\end{problem}

The Lazy Paver problem is known to be NP-complete
(\cite{EricksonRuskey2013}).  That is, \probref{lazy} is as difficult
to solve, in terms of the number of corners in the input region, as
the most difficult decision problems whose answers can be checked in
polynomial time.

Once again, we set up the game as a covering problem, rather than a
decision problem, by placing an overlay on top of the Tatami Maker
modules.  The overlay has a rectilinear region cut out of it, as well
as the instruction to cover the region with dominoes.  A prototypical
version and instructions appear in \fref{lazyPaver}.

\subsection*{The Paving Consultant}
\label{sec:consultant}

Dr. Matt DeVos communicated this final decision problem to me, and it
is also playable as a covering game (private communication).

\begin{problem}
  \label{prob:consultant}
  Determine whether a partial covering of a given rectilinear region
  can be completed.
\end{problem}

A Paver, perhaps a lazy one, has abandoned her job, and left a
partially paved region.  The Paving Consultant must decide whether the
covering can be completed.  The puzzle can be set up with overlays, as
above, bearing a diagram of the partial covering and indications on
how to arrange the Tatami Maker Modules.  Tiles of a different colour
can be provided so that the player does not mix these up with the
tiles that are part of the solution (see \fref{consultant}).

\begin{figure}[h!t]
  \centering
  \subfigure[ ]{%
    \includegraphics[width=0.35\textwidth]{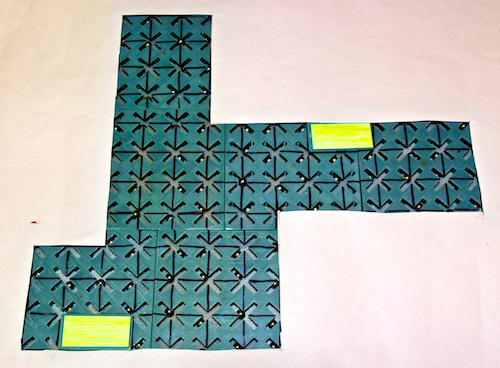}%
    \label{fig:consultantsetup}}%
  \hspace{0.2in}%
  \subfigure[ ]{%
    \includegraphics[width=0.35\textwidth]{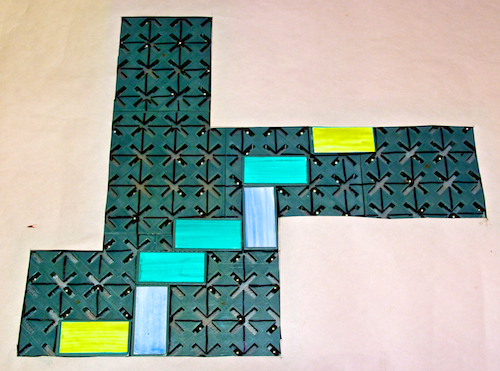}%
    \label{fig:consultantsolved}}%
  \caption{\emph{(a)}, A NO instance of the Paving Consultant.
    Placing more tiles, \emph{(b)}, hints at why this is so, and
    suggests an alternate definition for the tatami restriction.}
  \label{fig:consultant}
\end{figure}

\subsection*{Noku}
Noku is an adversarial game, perhaps antithetical to Oku, proposed by
Frank Ruskey (private communication).  Players alternately place tiles
onto the game board in order to win by forcing a position where their
opponent has no available move.  A computer analysis of some of Noku's
smallest game trees shows that Player $2$ can force a win on the
$2\times 6$ grid, which has $431949$ nodes, but the first player wins
for other similarly small cases, including $4\times 3$ and $4\times
4$.

\section*{Future}
I hope to produce more polished prototypes and test the playability of
these puzzles thoroughly.  It is unclear whether Tatami Maker is
suited to becoming a larger art installation, or a boxed puzzle set,
but it has already given tangible life to an otherwise abstract
research topic.



\end{document}